\documentclass[12pt]{article}
\usepackage{amsmath}
\newcommand{\be}{\begin{equation}}
\newcommand{\ee}{\end{equation}}
\newcommand{\ba}{\begin{eqnarray}}
\newcommand{\ea}{\end{eqnarray}}
\newcommand{\ban}{\begin{eqnarray*}}
\newcommand{\ean}{\end{eqnarray*}}

 \newcommand{\qed}{\hspace*{\fill}\rule{3mm}{3mm}\quad}
\newcommand{\Pf}{\noindent  {\em Proof.} }

\newcommand{\sect}[1]{\section{#1}  \setcounter{equation}{0}}

\newtheorem{lem}{Lemma}[section]

\begin{document}
\newtheorem{defn}[lem]{Definition}
\newtheorem{theo}[lem]{Theorem}
\newtheorem{cor}[lem]{Corollary}
\newtheorem{prop}[lem]{Proposition}
\newtheorem{rk}[lem]{Remark}
\newtheorem{ex}[lem]{Example}
\newtheorem{note}[lem]{Note}
\newtheorem{conj}[lem]{Conjecture}

\title{The Log Entropy Functional Along The Ricci Flow}
\author{Rugang Ye \\ {\small Department  of Mathematics} \\
{\small University of California, Santa  Barbara}}
\date{August 5, 2007}
\maketitle

\noindent 1. Introduction \\
2. The log entrop functional \\
3. Monotonicity of the log entropy and the logarithmic Sobolev constant \\
4. Proofs \\

\sect{Introduction}

In [P], Perelman introduced the entropy functional ${\mathcal W}(g,f,\tau)$ and established its monotonicity
along the Ricci flow.
An important application of Perelman's entropy monotonicity is for establishing
the $\kappa$-noncollpasing property of the Ricci flow under a finite upper bound for the
time, see [P] and [Y1].  In [Y2],[Y3] and [Y4], logarithmic Sobolev inequalities along the Ricci flow in
all dimensions $n \ge 2$
were obtained using Perelman's entropy monotonicity, which lead to Sobolev inequalities and $\kappa$-noncollpasing estimates.
In particular, uniform estimates without any condition on the time were obtained
 under the assumption that the first eigenvalue of the operator
$-\Delta+\frac{R}{4}$ is nonnegative at the start of the Ricci flow.

In this paper we introduce a new entropy functional which we call the {\it log entropy} because
of the appearance of an additional logarithmic operation compared with Perelman's
entropy functional.  Unlike Perelman's entropy functional, the parameter
$\tau$ does not appear in the log entropy. (A term involving the time $t$ enters into
the formula of the adjusted log entropy. However, its role is very different from
that played by
$\tau$ in Perelman's entropy.)
  On the other hand, the form of the
log entropy is directly related to the log gradient version of the logarithmic
Sobolev inequality, see [Y2, Appendix A].  Based on Perelman's entropy monotonicity we'll
establish the monotonicity of the log entropy (or the adjusted log entropy).  Because of the close relation between the
log entropy and the log gradient version of the logarithmic Sobolev inequality, we
can then show that the log gradient version  of the logarithmic Sobolev
inequality improves along the Ricci flow.  (This also leads to an improvement of the logarithmic Sobolev inequalities
along the Ricci flow obtained in [Y2], [Y3] and [Y4].) The Sobolev inequalities also improve along the  Ricci flow in an indirect way because they follow from the logarithmic Sobolev inequality.

\sect{The log entropy functional}

Let $M$ be a closed manifold of dimension $n\ge 1$.
\\

\noindent {\bf Definition 1}  We define the {\it log entropy} functional as follows

\ba \label{log1}
{\mathcal Y}_0(g, u)=-\int_M u^2 \ln u^2 dvol +\frac{n}{2} \ln  \left( \int_M (|\nabla u|^2+ \frac{R}{4}u^2)dvol
\right),
\ea
 where $g$ is a smooth metric on $M$ and $u \in W^{1,2}(M)$ satisfies
\ba \label{eigenpositive}
\int_M (|\nabla u|^2+ \frac{R}{4} u^2) dvol >0.
\ea
Here, all geometric quantities are associated with $g$.

More generally, we define the {\it log entropy} functional with {\it remainder} $a$ as follows
\ba
{\mathcal Y}_{a}(g, u)=-\int_M u^2 \ln u^2 dvol +\frac{n}{2} \ln  \left( \int_M (|\nabla u|^2+ \frac{R}{4}u^2)dvol
+a
\right).
\ea

Next we introduce the {\it adjusted log entropy} which  depends on an additional parameter $t$.
\\

\noindent {\bf Definition 2} We define the {\it adjusted log entropy} with remainder $a$ as follows
\ba
{\mathcal Y}_{a}(g, u, t)=-\int_M u^2 \ln u^2 dvol +\frac{n}{2} \ln  \left( \int_M (|\nabla u|^2+ \frac{R}{4}u^2)dvol
+a\right) +4at.
\ea
Obviously, ${\mathcal Y}_a(g, u)={\mathcal Y}_a(g, u, 0)$. \\

We'll need the following notation which is used in [Y2]. For a metric $g$ on $M$,
let $\lambda_0(g)$ denote  the first eigenvalue of the operator
$-\Delta+\frac{R}{4}$ for $g$. \\

\sect{Monotonicity of the log entropy and the logarithmic Sobolev constant}

Now we consider a smooth solution $g=g(t)$ of the Ricci flow on $M \times [\alpha, T)$ for some
$\alpha<T$, where $\alpha$ is finite.  Let $u=u(t)$ be a smooth positive solution of the backward evolution equation
\ba \label{uequation}
\frac{\partial u}{\partial t}=-\Delta u-\frac{|\nabla u|^2}{u}+\frac{R}{2}u
\ea
such that the normalization condition
\ba \label{u-1}
\int_M u^2 dvol=1
\ea
holds true for all $t$.
Equivalently, $v=v(t) \equiv u^2(t)$ is a smooth positive solution of the conjugate heat equation of Perelman
\ba
\frac{\partial v}{\partial t}=-\Delta v+R v
\ea
such that
\ba \label{v-1}
\int_M v dvol=1
\ea
for all $t$. Since
\ba
\frac{d}{dt} \int_M u^2 dvol=\frac{d}{dt} \int_M v dvol=0,
\ea
(\ref{u-1}) (or (\ref{v-1})) holds true for all $t$ iff it holds true for
any one value of $t$.

\begin{theo} \label{entropymonotone} Assume that $a >-\lambda_0(g(\alpha))$. Then
${\mathcal Y}_a(t) \equiv {\mathcal Y}_{a}(g(t), u(t), t)$ is nondecreasing. Indeed, we have
\ba \label{monotone1}
\frac{d}{dt}{\mathcal Y}_{a} &\ge& \frac{n}{4\omega}
\int_M |Ric-2 \frac{\nabla^2 u}{u}+2 \frac{\nabla u \otimes \nabla u}{u^2}
-\frac{4 \omega}{n} g|^2 u^2dvol \nonumber \\
&=&  \frac{n}{4\omega} \int_M |Ric+\nabla^2 f-\frac{4\omega}{n} g|^2 e^{-f} dvol ,
\ea
where $u=e^{-f/2}$ and
\ba
\omega=\omega(t)=a+\int_M (|\nabla u|^2+\frac{R}{4}u^2)dvol|_{t},
\ea
which is positive. 
\end{theo}

Note that if ${\mathcal Y}(t_2)={\mathcal Y}(t_1)$ for some $t_2>t_1$, then
the above monotonicity inequality (or (\ref{integral})) implies that 
\ba
Ric+\nabla^2 f-\frac{1}{2(t_2-t+\sigma)} g=0
\ea
on $[t_1, t_2]$ (and hence on $[t_1, T]$), i.e. $g$ is a gradient shrinking
soliton, where $\sigma=\frac{n}{8\omega(t_2)}$. (Note that the $f$ here is different from the $f$ employed in the proof of Theorem \ref{entropymonotone}
given below. This $f$ is used for the purpose of simplifying the expressions in the above formulas.) 

Next we define for each $a>-\lambda_0(g)$ the {\it logarithmic Sobolev constant with the
$a$-adjusted scalar curvature
potential}
\ba
C_{S, log, a}(M, g)&=&\inf\{-\int_M u^2 \ln u^2 dvol +\frac{n}{2} \ln  \left( \int_M (|\nabla u|^2+ (\frac{R}{4}+a)u^2)dvol
\right):  \nonumber \\
&& u \in W^{1,2}(M), \int_M u^2 dvol=1 \}.
\ea

In other words, $C_{S, log, a}(M,g)$ is the optimal constant (i.e. the maximal possible constant) such that the logarithmic
Sobolev inequality
\ba
\int_M u^2 \ln u^2 dvol \le \frac{n}{2} \ln  \left( \int_M (|\nabla u|^2+ \frac{R}{4}u^2)dvol+a\right)
-C_{S, log, a}(M,g)
\ea
holds true for all $u \in W^{1,2}(M)$ with $\int_M u^2 dvol=1$.  A natural question is how to estimate
$C_{S, log, a}(M,g )$ for a given $(M ,g)$.  The next proposition provides
such an estimate in terms of the modified Sobolev constant
$\tilde C_S(M,g)$, whose definition can be found  in [Y2].  A similar estimate for $n=2$ can be proved
in the same way by using the results in [Y3].   These estimates can be applied to
deal with the initial metric
when we study the Ricci flow.

\begin{prop} \label{logconstant} Assume $n \ge 3$. Let $g$ be a metric on $M$  and $a>-\lambda_0(g)$. If $a\ge -\frac{\min R^-}{4}
+\tilde C_S(M,g)^{-2} vol_g(M)^{-2/n}$, there holds
\ba
\int_M u^2 \ln u^2 dvol &\le& \frac{n}{2}\ln \left(\int_M (|\nabla u|^2+\frac{R}{4}u^2)dvol+a\right)
+\frac{n}{2} \ln (2 \tilde C_S(M,g)^2) \nonumber \\
\ea
for all $u \in W^{1,2}(M)$ with $\int_M u^2 dvol=1$.  If $a< -\frac{\min R^-}{4}
+\tilde C_S(M,g)^{-2} vol_g(M)^{-2/n}$, there holds
\ba
\int_M u^2 \ln u^2 dvol &\le& \frac{n}{2} \ln \left(\int_M (|\nabla u|^2+\frac{R}{4}u^2)dvol+a\right)
+\frac{n}{2}\ln (1+\frac{B}{a+\lambda_0(g)})\nonumber \\
&&+\frac{n}{2} \ln (2 \tilde C_S(M,g)^2),
\ea
for all $u \in W^{1,2}(M)$ with $\int_M u^2 dvol=1$, where
\ba
B=-a-\frac{\min R^-}{4}
+\frac{1}{\tilde C_S(M,g)^2 vol_g(M)^{2/n}}.
\ea
In other words, we have $C_{S, log, a}(M, g) \ge -\frac{n}{2} \ln (2 \tilde C_S(M,g))$ in the former case
and $C_{S, log, a}(M, g) \ge -\frac{n}{2}\ln (1+\frac{B}{a+\lambda_0(g)})-\frac{n}{2} \ln (2 \tilde C_S(M,g)^2)$
in the latter case.
\end{prop}
\Pf
We have by [Y2, Theorem 3.3]
\ba
\int_M u^2 \ln u^2 dvol &\le& \frac{n}{2} \ln \left(\tilde C_S(M,g) \|\nabla u\|_2+ \frac{1}{vol_g(M)^{1/n}}\right)^2
\nonumber \\
&\le&  \frac{n}{2} \ln (2\tilde C_S(M,g)^2) + \frac{n}{2} \ln \left(\int_M |\nabla u|^2 dvol+ \frac{1}{C_S(M,g)^2 vol_g(M)^{2/n}}\right)
\nonumber \\
&\le&  \frac{n}{2}\ln \left(\left[\int_M (|\nabla u|^2+\frac{R}{4}u^2)dvol+a\right]
+B\right)+\frac{n}{2} \ln (2\tilde C_S(M,g)^2), \nonumber \\
\ea
where
\ba
B=-a-\frac{\min R^-}{4}
+\frac{1}{\tilde C_S(M,g)^2 vol_g(M)^{2/n}}.
\ea
Hence we deduce in the case $B\le 0$, i.e. $a\ge -\frac{\min R^-}{4}
+\tilde C_S(M,g)^{-2} vol_g(M)^{-2/n}$
\ba
\int_M u^2 \ln u^2 dvol &\le& \frac{n}{2}\ln \left(\int_M (|\nabla u|^2+\frac{R}{4}u^2)dvol+a\right)
+\frac{n}{2} \ln (2 \tilde C_S(M, g)^2). \nonumber \\
\ea
In the case $B>0$, i.e. $a< -\frac{\min R^-}{4}
+\tilde C_S(M,g)^{-2} vol_g(M)^{-2/n}$, we consider the function $y=\ln(x+B)-\ln x$ for
$x \ge a+\lambda_0(g)$.  Clearly, $y$ is decreasing. Hence its maximum is
$y(a+\lambda_0(g))=\ln (a+\lambda_0(g)+B)-\ln (a+\lambda_0(g)).$  It follows that
$\ln (x+B) \le \ln x +\ln (a+\lambda_0(g)+B)-\ln (a+\lambda_0(g))$ for all $x \ge a+\lambda_0(g)$.
Consequently, there holds
\ba
\int_M u^2 \ln u^2 dvol &\le& \frac{n}{2} \ln \left(\int_M (|\nabla u|^2+\frac{R}{4}u^2)dvol+a\right)
+\frac{n}{2}\ln (1+\frac{B}{a+\lambda_0(g)})\nonumber \\
&&+\frac{n}{2} \ln (2\tilde C_S(M,g)^2).
\ea

The monotonicity of the log entropy functional given in Theorem \ref{entropymonotone} leads to the
monotonicity of the log Sobolev constant (or the adjusted log Sobolev constant), as formulated in the next theorem.

\begin{theo}  \label{logconstantmonotone} The adjusted logarithmic Sobolev inequality improves along the Ricci flow.
More precisely, $C_{S, log, a}(M, g(t))+4at$ is nondecreasing along an arbitrary
smooth solution $g(t)$ of the Ricci flow on $M$, provided that $a$ is greater than
the negative multiple of the $\lambda_0$ of the initial metric.  In particular, the logarithmic Sobolev
inequality improves along the Ricci flow, i.e.
$C_{S, log, 0}(M, g(t))$ is nondecreasing along the Ricci flow, provided that
$\lambda_0>0$ at the start.
\end{theo}

Next we handle the case $\lambda_0=0$ at the start.

\begin{theo} \label{zerocase} Assume $\lambda_0=0$ at the start.   Then either $g=g(t)$ is a gradient soliton, in which case
the logarithmic Sobolev constant $C_{S, log, a}(M, g(t))$ is independent of $t$ for any given $a>0$;
or $C_{S, log, 0}(M, g(t))$ is nondecreasing on $[\epsilon, T)$ for each $\epsilon>0$.
\end{theo}

Combining Theorem \ref{logconstantmonotone} and Theorem \ref{zerocase} with Proposition \ref{logconstant}
(and the corresponding result in dimension 2)  we then
arrive at a log gradient version
of the logarithmic Sobolev inequality along the Ricci flow, which improves the corresponding
results in [Y2], [Y3] and [Y4].

\sect{The proofs}

Now we proceed to prove Theorem \ref{entropymonotone}, Theorem \ref{logconstantmonotone} and
Theorem \ref{zerocase}.
We consider Perelman's entropy functional
\ba
{\mathcal W}(g, f, \tau)=
\int_M \left[ \tau(R+|\nabla f|^2)+f-n\right] \frac{e^{-f}}{(4\pi\tau)^{\frac{n}{2}}} dvol.
\ea
Setting
\ba
u=\frac{e^{-f/2}}{(4\pi \tau)^{n/4}},
\ea
i.e.
\ba
f=-\ln u^2-\frac{n}{2} \ln \tau-\frac{n}{2} \ln (4 \pi),
\ea
we have
\ba \label{entropy}
{\mathcal W}(g, f, \tau)&=& -\int_M u^2 \ln u^2 dvol+\tau \int_M (4|\nabla u|^2+R u^2)dvol-\frac{n}{2}\ln \tau \nonumber \\
&& -\frac{n}{2}\ln(4\pi)-n \nonumber \\
&=& -\int_M u^2 \ln u^2 dvol+\tau \left[\int_M (4|\nabla u|^2+R u^2)dvol+4a\right]-\frac{n}{2}\ln \tau \nonumber \\
&& -\frac{n}{2}\ln(4\pi)-n -4a \tau \nonumber \\
&=& -\int_M u^2 \ln u^2 dvol+(4\tau) \left[\int_M (|\nabla u|^2+\frac{R}{4} u^2) dvol+a\right]-\frac{n}{2}\ln (4\tau) \nonumber \\
&& -\frac{n}{2}\ln \pi-n -a (4\tau)
\ea
for an arbitary constant $a$.

\begin{lem} \label{minimumlemma1} Let a metric $g$ on $M$ be given. Assume $a>-\lambda_0(g)$. Let
$u \in W^{1,2}(M)$ with $\int_M u^2 dvol=1$.  Then the minimum of the function
\ba
h(s)=s \left[\int_M (|\nabla u|^2+\frac{R}{4} u^2)dvol+a\right]-\frac{n}{2}\ln s
\ea
for $s>0$
is given by
\ba
\min h=\frac{n}{2}\ln \left(\int_M (|\nabla u|^2+\frac{R}{4} u^2)dvol+a\right) +\frac{n}{2}(1-
\ln \frac{n}{2})
\ea
and is achieved  at the unique minimum point
\ba
s=\frac{n}{2} \left(\int_M (|\nabla u|^2+\frac{R}{4} u^2)dvol+a\right)^{-1}.
\ea
\end{lem}
\Pf Set $A=\int_M (|\nabla u|^2+\frac{R}{4} u^2)dvol+a$. Then $h(s)=As-
\frac{n}{2}\ln s$.  Since $a>-\lambda_0(g)$, we have
$A>0$.  Hence we have $h(s) \rightarrow \infty$ as $s \rightarrow \infty$.
We also have $h(s) \rightarrow \infty$ as $s \rightarrow 0$. Consequently, $h$ achieves its minimum
somewhere. Since $h'(s)=A-\frac{n}{2s}$, the minimum is achieved at the unique
minimum point $s=
\frac{n}{2A}$, and then the minimum is $h(\frac{n}{2A})=\frac{n}{2}\ln A+\frac{n}{2}(1-
\ln \frac{n}{2})$.  \qed \\

The following lemma is a simple corollary of the above lemma.

\begin{lem} \label{minimumlemma2}  Assume $a>-\lambda_0(g)$. Then there holds
for each $\tau>0$
\ba
{\mathcal W}(g, f, \tau) \ge
-\int_M u^2 \ln u^2 dvol+ \frac{n}{2}\ln \left(\int_M (|\nabla u|^2+\frac{R}{4} u^2)dvol+a\right)
-4a\tau +b(n), \nonumber \\
\ea
where
\ba
b(n)= -\frac{n}{2}\ln \pi-\frac{n}{2}(1+
\ln \frac{n}{2}).
\ea
Moreover, we have
\ba
{\mathcal W}(g, f, \frac{n}{8}\omega(g,u,a)^{-1})
&=&-\int_M u^2 \ln u^2 dvol+ \frac{n}{2}\ln \left(\int_M (|\nabla u|^2+\frac{R}{4} u^2)dvol+a\right)
\nonumber \\
&&-\frac{na}{2}\omega(g,u,a)^{-1} +b(n),
\ea
where
\ba
\omega(g,u,a)=\int_M (|\nabla u|^2+\frac{R}{4} u^2)dvol+a
\ea

\end{lem}

Now let $g=g(t)$ be a smooth solution of the
Ricci flow on $M \times [\alpha, T)$ and $u=u(t)$ as postulated for Theorem \ref{entropymonotone}. Fix
$\alpha \le t_1<t_2 <T$ and define for a given $\sigma>0$
\ba
\tau=\tau(t)=t_2-t+\sigma.
\ea
We define $f=f(t)$ by
\ba
f=-\ln u^2-\frac{n}{2} \ln \tau-\frac{n}{2} \ln (4 \pi),
\ea
i.e.
\ba
u=\frac{e^{-f/2}}{(4\pi \tau)^{n/4}}.
\ea
Then $f$ satisfies the backward evolution equation
\ba
\frac{\partial f}{\partial t}=-\Delta f+|\nabla f|^2-R+\frac{n}{2\tau}.
\ea

By Perelman's entropy monotonicity formula we have  for $g=g(t), f=f(t)$ and $\tau=\tau(t)$
\ba \label{monotone2}
\frac{d}{dt} {\mathcal W}(g,f,\tau)=2\tau \int_M |Ric+\nabla^2 f-\frac{1}{2\tau} g|^2 \frac{e^{-f}}{(4\pi\tau)^{\frac{n}{2}}} dvol
\ea
on $[t_1, t_2]$. It follows that
\ba \label{monotone3}
{\mathcal W}(g(t_2), f(t_2), \sigma) &=& {\mathcal W}(g(t_1), f(t_1), t_2-t_1+\sigma)
\nonumber \\
&& + 2 \int_{t_1}^{t_2} \tau \int_M |Ric+\nabla^2 f-\frac{1}{2\tau} g|^2 \frac{e^{-f}}{(4\pi\tau)^{\frac{n}{2}}} dvol dt.
\ea

\noindent {\bf Proof of Theorem \ref{entropymonotone} }\\

Choosing $\sigma=\frac{n}{8} \omega(g(t_2), u(t_2), a)^{-1}$ in (\ref{monotone3}) we deduce
\ba
&& -\int_M u^2 \ln u^2 dvol|_{t_2}+\frac{n}{2}\ln \left(\int_M (|\nabla u|^2+\frac{R}{4} u^2)dvol+a\right)|_{t_2}
-4a\sigma+b(n)
\nonumber \\ && = {\mathcal W}(g(t_1), f(t_1), t_2-t_1+\sigma)
\nonumber \\
&& + 2 \int_{t_1}^{t_2} (t_2-t+\sigma) \int_M |Ric+\nabla^2 f-\frac{1}{2(t_2-t+\sigma)} g|^2 \frac{e^{-f}}{(4\pi(t_2-t+\sigma))^{\frac{n}{2}}} dvol dt
\nonumber \\
&& \ge -\int_M u^2 \ln u^2 dvol |_{t_1}+ \frac{n}{2}\ln \left(\int_M (|\nabla u|^2+\frac{R}{4} u^2)dvol+a\right)|_{t_1}
\nonumber \\
&& -4a(t_2-t_1+\sigma)+b(n) \nonumber \\
&&   +2 \int_{t_1}^{t_2} (t_2-t+\sigma) \int_M |Ric+\nabla^2 f-\frac{1}{2(t_2-t+\sigma)} g|^2 \frac{e^{-f}}{(4\pi(t_2-t+\sigma))^{\frac{n}{2}}} dvol dt
\nonumber \\
\ea

It follows that
\ba \label{integral}
&&{\mathcal Y}_a(g(t_2), u(t_2), t_2) \ge {\mathcal Y}_a(g(t_1), u(t_1), t_1)
\nonumber \\
&& + 2 \int_{t_1}^{t_2} (t_2-t+\sigma) \int_M |Ric+\nabla^2 f-\frac{1}{2(t_2-t+\sigma)} g|^2 \frac{e^{-f}}{(4\pi(t_2-t+\sigma))^{\frac{n}{2}}} dvol dt,\nonumber \\
\ea
which leads to
\ba
\frac{d}{dt} {\mathcal Y}_a(g(t), u(t), t) &\ge& 2 \sigma \int_M |Ric+\nabla^2 f-\frac{1}{2\sigma} g|^2 \frac{e^{-f}}{(4\pi\sigma)^{\frac{n}{2}}} dvol \nonumber \\
&=& \frac{n}{4\omega}
\int_M |Ric-2 \frac{\nabla^2 u}{u}+2 \frac{\nabla u \otimes \nabla u}{u^2}
-\frac{4 \omega}{n} g|^2 u^2dvol.\nonumber
\ea
This implies that ${\mathcal{Y}}_a(g(t), u(t), t)$ is nondecreasing.  \qed \\

\noindent {\bf Remark}  With some more work, we can
show that ${\mathcal{Y}}_a(g(t), u(t), t)$ is actually strictly increasing unless $a=0$.  \\

\noindent {\bf Proof of Theorem \ref{logconstantmonotone} } Consider $t_1<t_2$.
Given $\epsilon>0$, we choose a positive smooth function $u_2$ such that
$\int_M u_2^2 dvol=1$ at $t=t_2$ and ${\mathcal{Y}}_a(g(t_2), u_2) \le
C_{S, log, a}(M, g(t_2))+\epsilon$. Then we let $u=u(t)$ be the positive smooth solution of the
equation (\ref{uequation}) on $[t_1, t_2]$ with $u(t_2)=u_2$. By Theorem \ref{entropymonotone}
we then have
\ba
C_{S,log,a}(M, g(t_2)) + \epsilon +4at_2 &\ge& {\mathcal Y}_a(g(t_2), u(t_2)) +4at_2 \ge {\mathcal Y}_a(g(t_1), u(t_1))+4at_1 \nonumber \\
&\ge&
C_{S, log, a}(M ,g(t_1))+4at_1.
\ea
Since $\epsilon$ is arbitary, we arrive at
\ba
C_{S,log,a}(M, g(t_2))+4at_2 \ge C_{S, log, a}(M ,g(t_1))+4at_1.
\ea
\qed \\

\noindent {\bf Proof of Theorem \ref{zerocase} } This follows from the arguments
in [Y4] and Theorem \ref{logconstantmonotone}. \qed \\

\end{document}